\pdfoutput=1
\documentclass[bezier,amstex, oneside,reqno, twocolumn]{amsproc}
\usepackage{amsmath, amssymb, amsthm, verbatim, euscript, amscd, textcomp, mathrsfs, caption}
\usepackage[dvips]{graphicx}
\usepackage{epsfig}
 \usepackage[all,cmtip]{xy}
\usepackage{color}

\numberwithin{figure}{section}

\def\ie{\emph{i.e., }}
\def\eg{\emph{e.g., }}

\def\R{\mathbb R}
\def\Z{\mathbb Z}

\def\T{\mathbb T}

\def\eps{\varepsilon}

\newtheorem{theorem}{Theorem}[section]
\newtheorem{conj}[theorem]{Conjecture}
\newtheorem{problem}[theorem]{Problem}

 \theoremstyle{remark}
\newtheorem{remark}[theorem]{Remark}

\theoremstyle{remark}

\begin{document}
\author{Andrey Gogolev, Itai Maimon and  Aleksey N. Kolmogorov}
\title[Gibbs $u$-measures on $\T^3$]{A numerical study of Gibbs $u$-measures for partially hyperbolic diffeomorphisms on $\T^3$.}

\begin{abstract}
We consider a hyperbolic automorphism $A\colon\T^3\to\T^3$ of the 3-torus whose 2-dimensional unstable distribution splits into weak and strong unstable subbundles. We unfold $A$ into two one-parameter families of Anosov diffeomorphisms --- a conservative family and a dissipative one. For diffeomorphisms in these families we numerically calculate the strong unstable manifold of the fixed point. Our calculations strongly suggest that the strong unstable manifold is dense in $\T^3$. Further, we calculate push-forwards of the Lebesgue measure on a local strong unstable manifold. These numeric data indicate that the sequence of push-forwards converges to the SRB measure.
\end{abstract}
\date{}
 \maketitle

\section{Introduction}
\subsection{The setting}
Consider the 3-dimensional torus $\T^3=\R^3/\Z^3$ equipped with the standard $(x,y,z)$ coordinates and a hyperbolic automorphism $A\colon \T^3\to \T^3$ induced by the following integral matrix with determinant 1
$$
A=
\begin{pmatrix}
 2 & 1 & 0 \\
 1 & 2 & 1 \\
 0 & 1 & 1 \\
\end{pmatrix}
$$
The eigenvalues of $A$ are real and approximately equal to $0.20$, $1.55$ and $3.25$. We denote the largest eigenvalue by $\lambda$, $\lambda\approx 3.25$, and corresponding eigenvector by $v$, $Av=\lambda v$,
$$
v\approx
\begin{pmatrix}
  0.80 \\
 1.00  \\
 0.45  \\
\end{pmatrix}
$$
We will view $A$ as a partially hyperbolic diffeomorphism whose center distribution is expanding. Further, we unfold $A$ into two families of partially hyperbolic diffeomorphisms: a {\it dissipative family}
\begin{multline}
\label{D}
f_{D, \eps}(x,y,z)=\\(2x+y  +\eps \sin(2\pi x), x+2y+z, y+ z)
\end{multline}
 and a {\it conservative family}
 \begin{multline}
 \label{C}
f_{C, \eps}(x,y,z)=\\(2x+y +\eps \sin(2\pi x),  x+2y+z+ \eps \sin(2\pi x), y+ z)
\end{multline}
It is well-known that for small values of $\eps>0$ the diffeomorphisms $f_{*,\eps}$, (here $*=D,C$), remain Anosov, and also partially hyperbolic (with weakly expanding center distribution). Hence diffeomorphisms $f_{*,\eps}$ leave invariant a one-dimensional strongly expanding foliation $W^{uu}_{f}$ whose expansion rate is close to $\lambda$. Note that the point $p=(0,0,0)$ is fixed by all diffeomorphisms in the families. 

\subsection{Preview of the results and conjectures} 
\begin{figure}
\includegraphics[width=70mm]{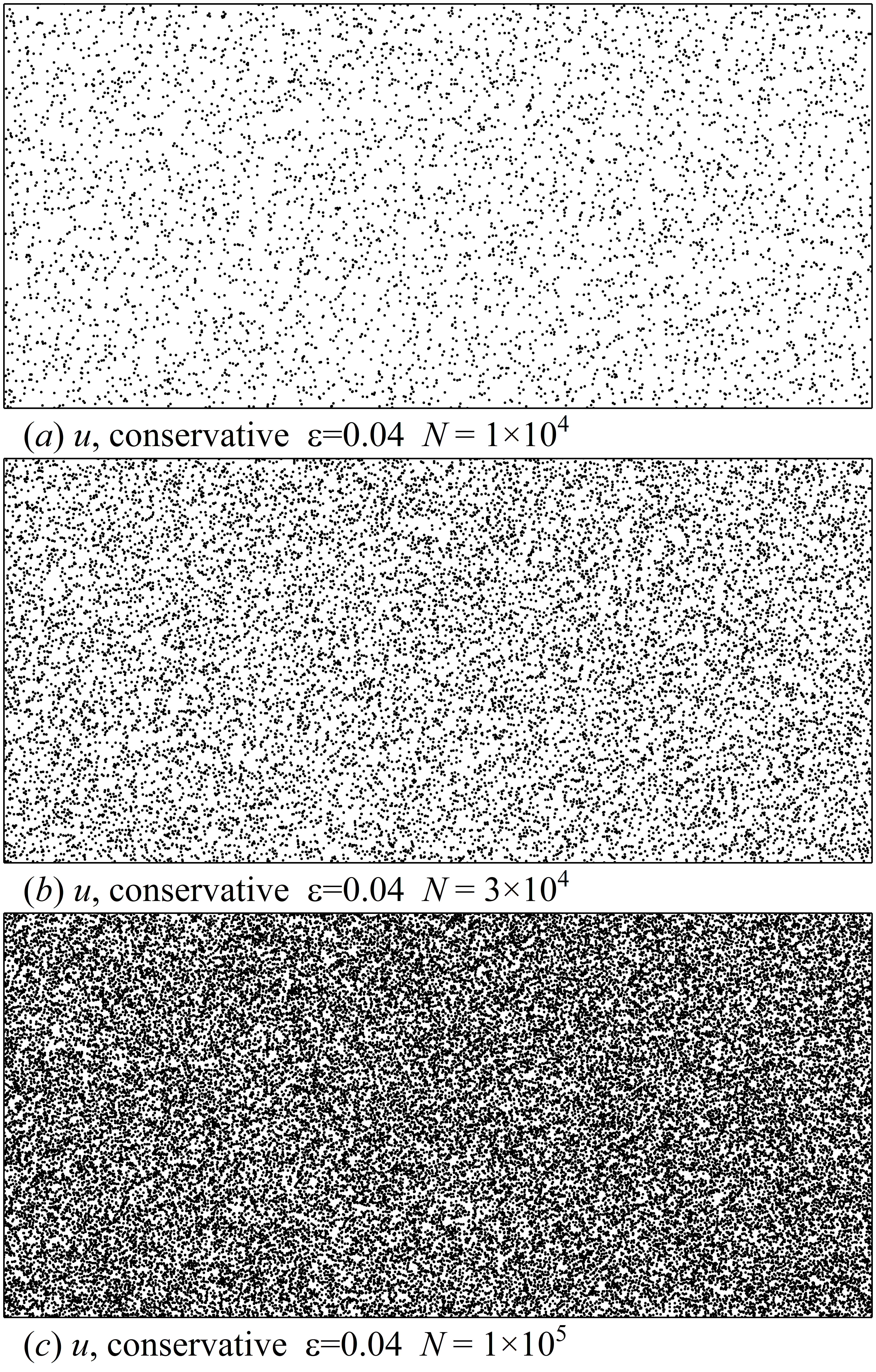}
\caption{Intersection points of the strong unstable manifold and $\T^2$ transversal for $f_C$ with $\eps = 0.04$. The snapshots are shown for the first $N = 10,000$, $30,000$ and $100,000$ points.}
\label{fig1}
\end{figure}

\begin{figure}
\includegraphics[width=70mm]{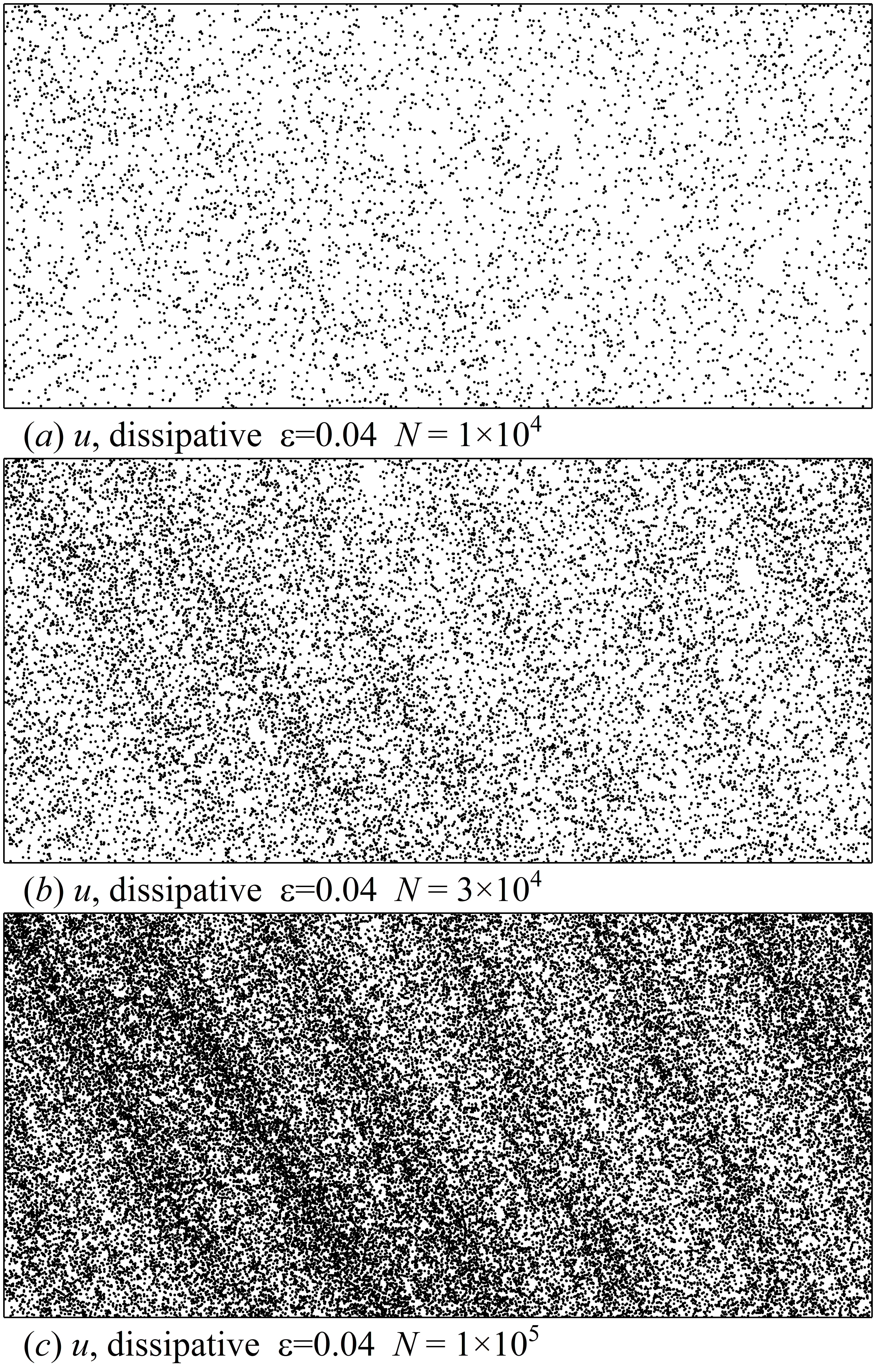}
\caption{Intersection points of the strong unstable manifold and $\T^2$ transversal for  $f_D$ with $\eps = 0.04$. The snapshots are shown for the first $N = 10,000$, $30,000$ and $100,000$ points.}
\label{fig2}
\end{figure}

We performed a very accurate (albeit non-rigorous) numerical calculations of the the finite-length strong unstable manifolds $W^{uu}_{*,\eps}(p, R)$ which pass through $p$, up to length $R\approx 1.3\cdot10^8$. These numerical calculations strongly support the following conjecture.
\begin{conj}
\label{conj1}
For all analytic diffeomorphisms $f$ in a sufficiently small neighborhood of $A$ the strong unstable foliation $W^{uu}_f$ is transitive, \ie it has a dense leaf.
\end{conj}
We actually expect the foliation $W^{uu}_f$ to be minimal. However we did not calculate strong unstable leaves through non-periodic points because it is a much harder task. 
Figures~\ref{fig1} and~\ref{fig2} give a preview of our numerics in support of the above conjecture. These panels display first $N$ intersection points of the strong unstable manifold $W^{uu}(p)$ with the 2-torus given by $y=0$. We have calculated $10^8$ intersection points. In the figures we only show up to $200,000$  points because of large file size and because past $10^6$ points one only sees Malevich's black square.

\captionsetup{width=70mm}

\begin{remark}
In fact we only show roughly $N/2$ points because we only display the points of intersection in the ``one-half torus" given by $0\le z\le 1/2$. We will display all our data on $[0, 1]\times [0,1/2]\subset [0,1]\times[0,1]\simeq \T^2$ unless specified otherwise. Note that throughout the paper we maintain the convention to indicate the total number of points $N$ in the captions to the figures. Hence, if readers counts the points on a figure then they would get approximately $N/2$ points.
Displaying half of the torus helps to reduce file size.
Also note that all diffeomorphisms $f_{*,\eps}$ commute with the involution $i\colon(x,y,z)\mapsto (-x,-y,-z)$. It follows that the measures which we are interested in are invariant under $i$ and all conditional measures on the $\T^2$ transversal are invariant under $(x,y)\mapsto (-x,-y)$. 
\end{remark}

In general, given a partially hyperbolic diffeomorphism $f\colon M\to M$, one reason to be interested in minimal sets of its strong unstable foliation $W^{uu}_f$ is that minimal invariant sets support Gibbs $u$-measures associated to $W^{uu}_f$ of Pesin and Sinai~\cite{PS}. Gibbs $u$-measures are of great interest in partially hyperbolic dynamics because they govern statistical properties of the dynamical system~\cite{Dol1, Dol2}. Of course, in our setting the dynamical system is a transitive Anosov diffeomorphism which admits a unique SRB-measure and, hence, statistical properties are very well understood. However, perturbations of linear partially hyperbolic automorphisms are nice model examples where $u$-measures are not fully understood. We elaborate on our motivation to carry out the numerical study at the end of the introduction.

\begin{figure}
\includegraphics[width=70mm]{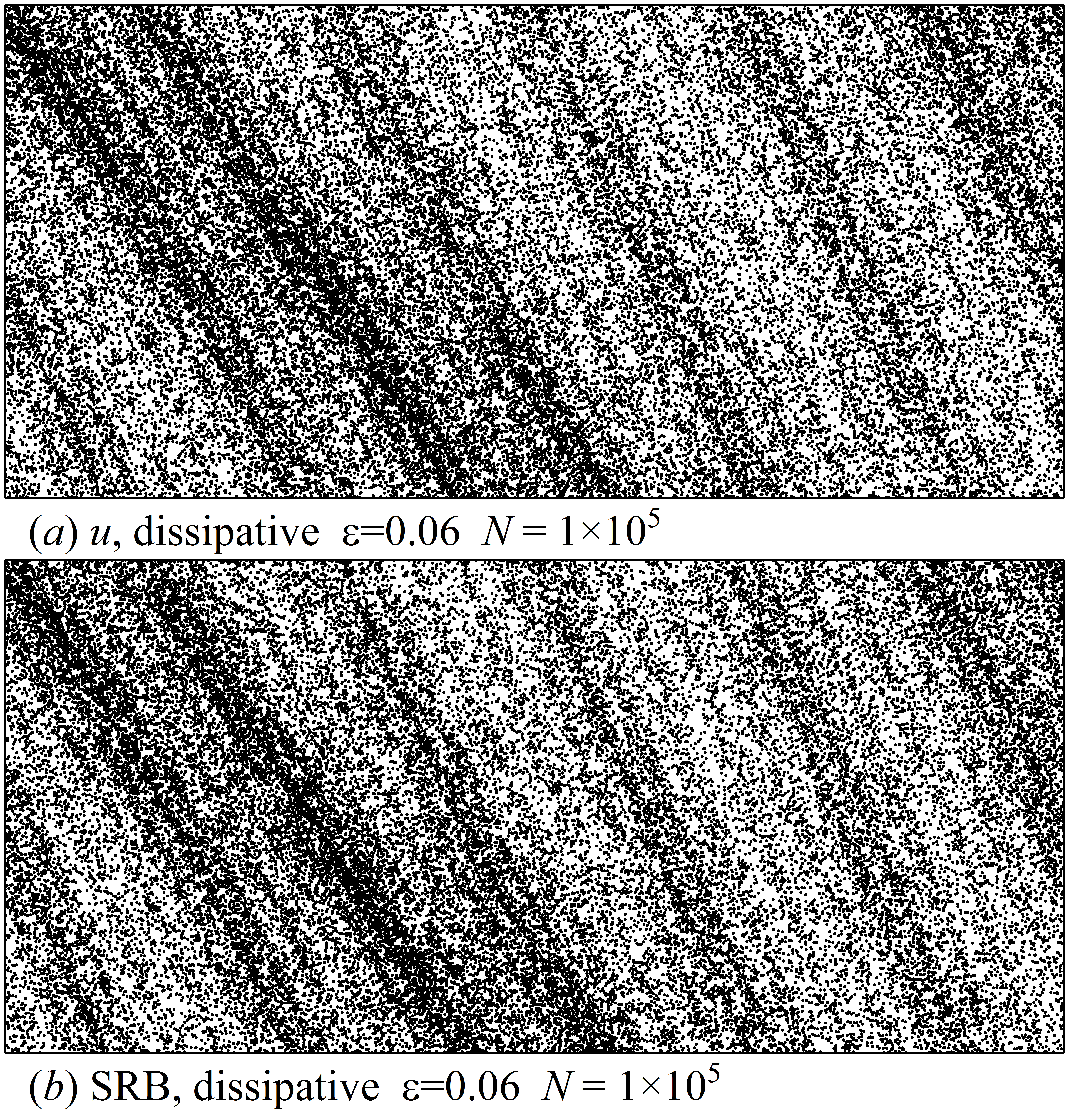}
\caption{Conditionals for $u$ and SRB measure of $f_D$ with $\eps = 0.06$.}
\label{fig3}
\end{figure}

We view diffeomorphisms given by~(\ref{D}) and~(\ref{C}) as partially hyperbolic diffeomorphisms with one-dimensional strong unstable subbundles. Recall that a {\it Gibbs $u$-measure} of a partially hyperbolic diffeomorphism $f\colon M\to M$ is an $f$-invariant measure $\mu$ whose conditional measures on strong unstable plaques are absolutely continuous with respect to the induced Riemannian volume on strong unstable plaques. Gibbs $u$-measures were introduced by Pesin and Sinai~\cite{PS}\footnote{Pesin and Sinai used a stronger definition which is equivalent to the one we give here, see~\cite[Chapter 11]{BDV}.} who also suggested a way to construct them as weak$^*$ partial limits of the sequence of averages 
\begin{equation}
\label{eq_averege}
\bar\nu^{uu}_K \stackrel{\mathrm{def}}{=} \frac{\nu^{uu}+f_*(\nu^{uu})+\ldots +f_*^{K-1}(\nu^{uu})}{K}, K\ge 1,
\end{equation}
where $\nu^{uu}$ is a singular measure (on $M$) given by induced Riemannian volume on a strong unstable plaque. In our setting we can take $\nu^{uu}$ to be the singular measure (on $M$) given by the Lebesgue measure on a small plaque of $W^{uu}(p)$ with one end point being $p$. Hence we amend our calculation of $W^{uu}(p)$ with a numeric calculation of the strong unstable Jacobians of $f^i$, $i\le K$, to obtain the averages numerically (more precisely, we look at the conditional measures of the averages on the 2-torus given by $y=0$). Even though our evidence is not completely conclusive (due to exponential error accumulation we do not go beyond $K=50$) we believe that the averages $\bar\nu_K^{uu}$ converge weakly. Further we calculate the SRB measure employing the zero-noise limit description of Young~\cite{Y2}.  The very different numeric procedures for calculating the $u$-measure and the SRB measure produce visually identical results for all values of $\eps$ as indicated on Figure~\ref{fig3}. Hence we cautiously conjecture the following.

\begin{conj}
\label{conj2}
For all analytic diffeomorphisms $f$ in a sufficiently small neighborhood of $A$ there exists a unique Gibbs $u$-measure (an $f$-invariant measure with absolutely continuos conditionals on strong unstable leaves) which then, of course, coincides with the SRB measure.
\end{conj}

This conjecture can be reformulated as follows: {\it for any analytic diffeomorphisms $f$ in a sufficiently small neighborhood of $A$ any $f$-invariant measure with absolutely continuous conditional measures on one-dimensional strong unstable plaques, in fact, has absolutely continuous conditional measures on two-dimensional unstable plaques.}

\subsection{Motivation}
\subsubsection{} Our initial interest in transitivity (or minimality) question of the strong unstable foliation came from work on smooth conjugacy of higher dimensional Anosov diffeomorphisms~\cite{Gog}. Transitivity of invariant expanding one-dimensional foliations (albeit not the strong unstable ones) played a key role in the arguments of~\cite{Gog}. Families of diffeomorphisms in dimension three which we consider in this paper is the simplest setting where transitivity (minimality) is not understood. 

We remark that minimality of the weak unstable foliation for $f_{*,\eps}$ follows easily from structural stability. Indeed the conjugacy between $A$ and $f$ sends the weak unstable foliation $W^{wu}_A$ to the weak unstable foliation $W^{wu}_f$ of $f$. It is well known (see, \eg~\cite{GG}) that in general the conjugacy does not respect the strong unstable foliation. In fact, according to~\cite{RGZ}, the conjugacy respects strong unstable foliations if and only if the strong unstable and stable distributions of $f$ integrate to an invariant foliation.

 Note also that here we consider an irreducible automorphism $A$ as a base-point for the families.
The situation  is quite different (but also poorly understood) in the case of reducible automorphisms. Indeed, for a reducible Anosov automorphism with weak-strong splitting in dimension 4 the closures of strong unstable leaves are 2-tori. And it is a very interesting question to investigate these closures and Gibbs $u$-measures for families which bifurcate into non-skew-product diffeomorphisms.

\subsubsection{} Minimality of strong unstable foliations of 3-dimensional partially hyperbolic diffeomorphisms was studied by Bonatti-D\'iaz-Ures~\cite{BDU}. In particular, results of~\cite{BDU} yield ``large" $C^1$-open sets of partially hyperbolic diffeomorphism with minimal strong stable and minimal strong unstable foliations. One of the crucial assumptions in the setting of~\cite{BDU} is existence of hyperbolic periodic points of different indices (that is, a index 1 and index 2). The main technique of~\cite{BDU} is construction of an invariant Morse-Smale section to the foliation. In our setting, when the center foliation is weakly expanding, this technique is not applicable. Hence our setting can be considered as a complementary one to the setting of~\cite{BDU}.

Recall that homotopy class of $A\colon\T^3\to\T^3$ contains the Ma\~n\'e's example. This is a robustly transitive diffeomorphism $f_M\colon\T^3\to\T^3$ which is partially hyperbolic but not Anosov~\cite{M}. To the best of our knowledge minimality of strong unstable foliation of $f_M$ is also an open problem. Thus, understanding strong unstable foliation of perturbations of $A$ and Ma\~n\'e's example is a prerequisite for the following problem, which is a special case of Problem~1.6 in~\cite{BDU}.
\begin{problem}
Consider the space of robustly transitive partially hyperbolic diffeomorphisms $f\colon \T^3\to\T^3$ which are homotopic to $A$, \ie the induced map $f_*$ on first homology group is given by $A$. Is strong unstable foliation $W^{uu}_f$ minimal? Or at least transitive? If not, then is minimality (transitivity) of $W^{uu}_f$ a $C^1$-open and dense property in this space?
\end{problem}

Related to this problem, Potrie asked whether transitivity (or chain-recurrence) follows from partial hyperbolicity of $f\colon \T^3\to\T^3$ in the homotopy class of $A$~\cite{P}. Further, Potrie proved that there exists a unique qausi-attractor for each such $f$. Note that the attractor must be saturated by leaves of $W^{uu}_f$. Hence, minimality of $W^{uu}_f$ would imply that the attractor is whole $\T^3$.\footnote{The first author believes that non-trivial trapping regions exist for some point-wise partially hyperbolic $f\colon\T^3\to\T^3$ in the homotopy class of $A$. If so, it would be interesting to investigate the structure of $W^{uu}_f$ and how the bifurcation happens.}

We also remark that robust minimality of strong unstable foliation was established in~\cite{PS} under so called SH-condition. This condition does not hold in our setting. Finally, minimal sets of strong unstable foliation can be analyzed better in $C^1$ generic setting, see~\cite[Section 5.3]{CP}.

\subsubsection{} It is interesting to understand the space of Gibbs $u$-measure $\textup{Gibbs}^u(f)$, its dependence on the diffeomorphism and what bifurcations can occur. Note that it is known that $\textup{Gibbs}^u(f)$ depends continuously on $f$ in $C^1$ topology~\cite{Yang} (see also~\cite[Chapter 9]{BDV}). Generalizing our numeric observation of uniqueness of the $u$ measure we ask the following question. 

\begin{problem} Let $f\colon M\to M$ be  a partially hyperbolic diffeomorphism of a 3-manifold $M$. Assume that $f$ admits two distinct ergodic $u$-measures $\mu_1$ and $\mu_2$. Is it true that $\textup{supp}(\mu_1)\neq\textup{supp}(\mu_2)$? That is, do they necessarily have distinct supports?
\footnote{It was suggested to us by Dmitry Dolgopyat that this question also makes sense in higher dimensions if one additionally assumes that $f$ is accessible (or considers $u$-measures supported on an accessibility class).}
\end{problem}

\subsection{Further discussion}

\begin{figure}
\includegraphics[width=70mm]{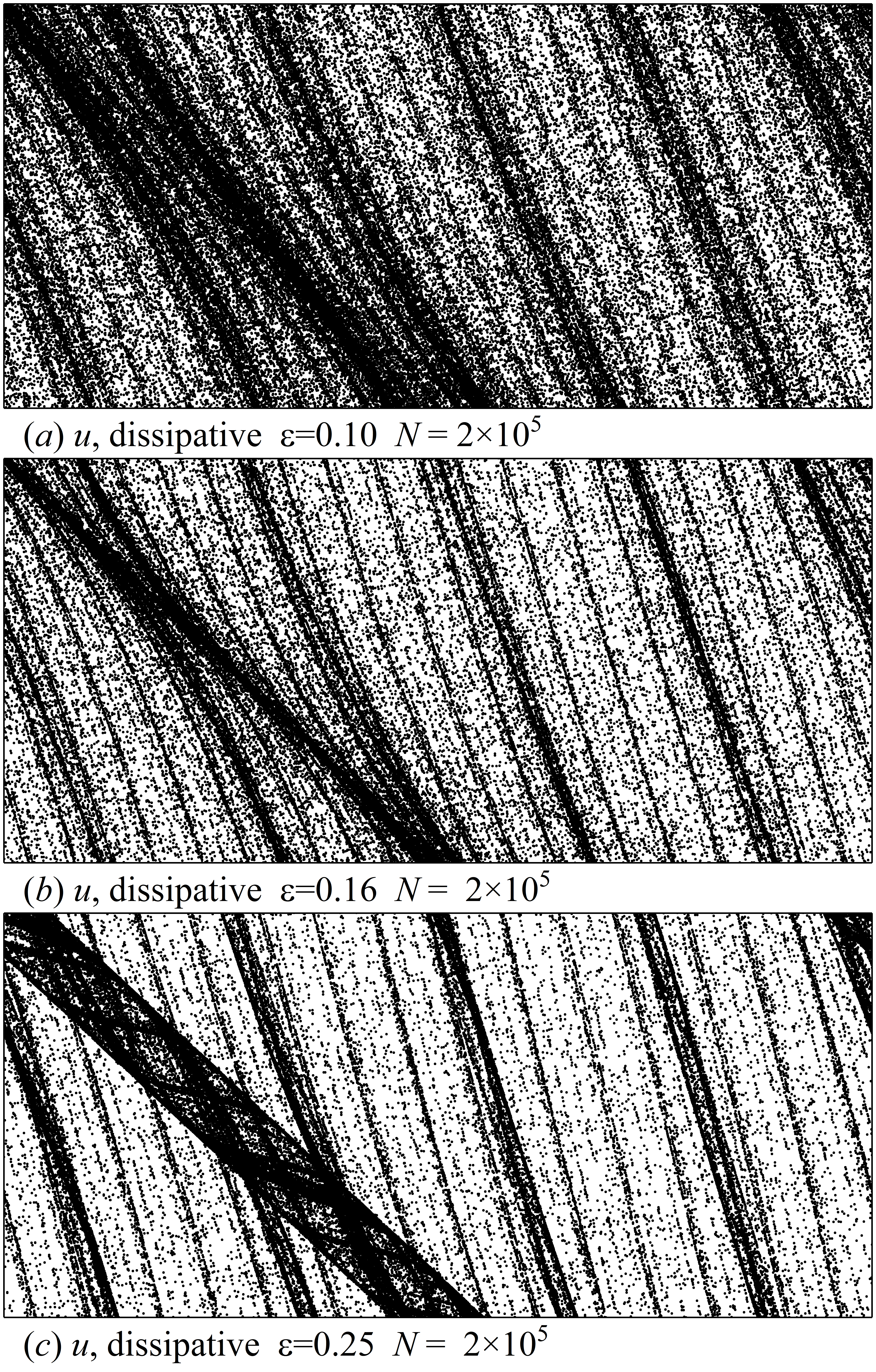}
\caption{Conditionals of the $u$-measure on $\T^2$ in the dissipative family $f_D$.}
\label{fig4}
\end{figure}

\subsubsection{} Our numerical evidence actually suggests that the push-forward measures $f^n_*\nu^{uu}$ converge to the SRB measure as $n\to\infty$. This was easier to detect than convergence of averages~(\ref{eq_averege}), which clearly converge slower. Note that convergence of $f^n_*\nu^{uu}$ to the unique Gibbs $u$-measure is a key assumption in the study of statistical properties of partially hyperbolic diffeomorphisms~\cite{Dol2}. This assumption is difficult to verify theoretically when dynamics along the center subbundle is non-linear.\footnote{However, for transitive Anosov diffeomorphisms mixing implies that $f^n_*\nu^{u}$ converges to the SRB measure, where $\nu^u$ is Lebesgue measure on an unstable plaque. This was explained to us by F. Rodriguez Hertz.}

\begin{figure}
\includegraphics[width=85mm]{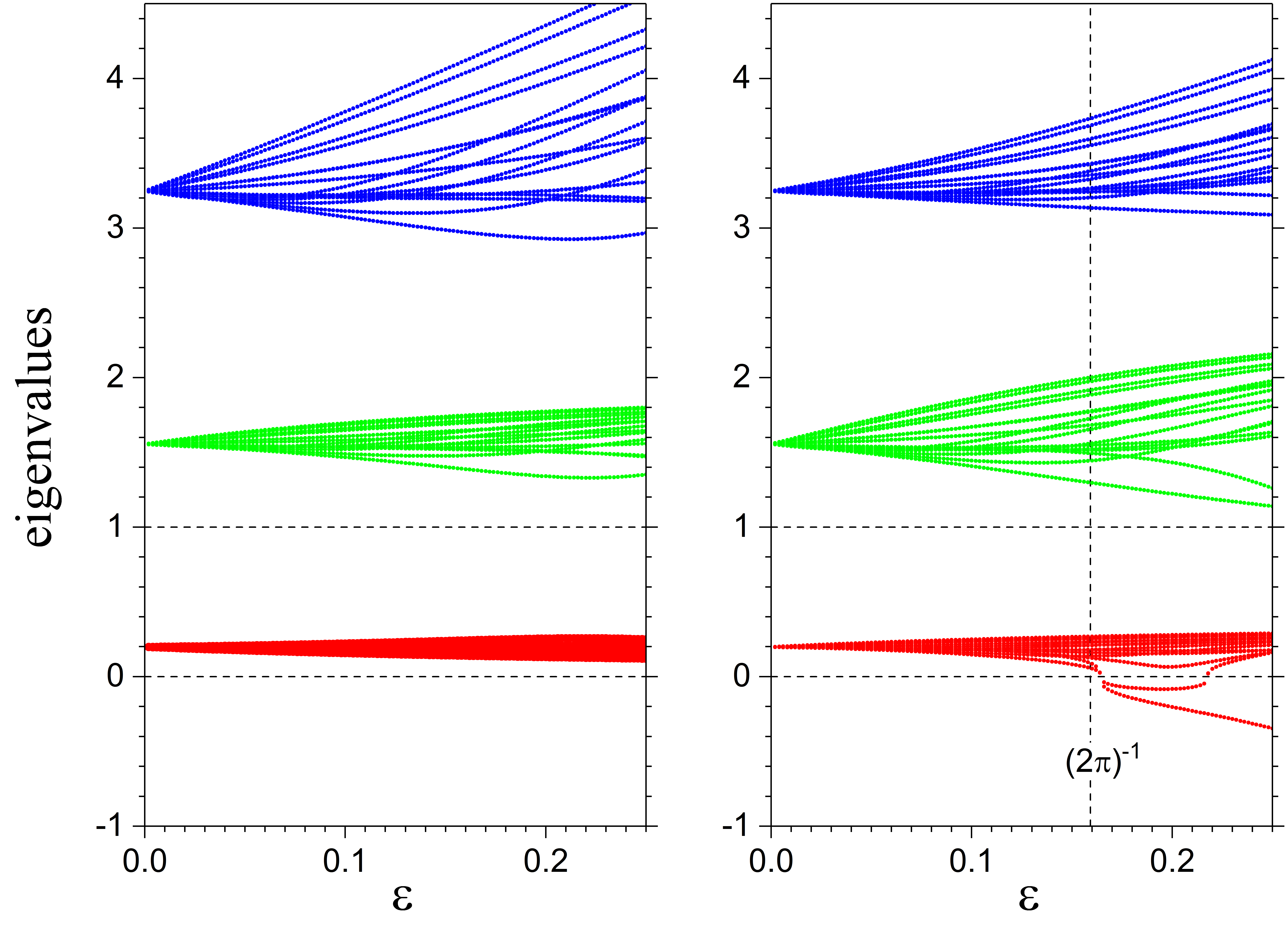}
\caption{Spectrum at the fixed point and $15$ orbits of period 3. Note that it is disjoint with horizontal line 1 and is confined to three disjoint bands, suggesting that diffeomorphisms $f_C$ and $f_D$ stay Anosov with weak-strong unstable splitting.}
\label{fig45}
\end{figure}

\subsubsection{} 

Another observation is that our numerical experiments suggest that the
$u$-measure coming from the strong unstable leaf through $p$ agrees
with the SRB measure well beyond the range of small $\eps$. (The
splitting at $p$ survives for all $\eps>0$.) This is indicated on
Figure~\ref{fig4}. At $\eps=\frac{1}{2\pi}\approx 0.159$ bifurcation
from diffeomorphisms to non-invertible maps occurs and the ``folding"
which happens beyond this parameter value is clearly visible on
Figure~\ref{fig4}. Note that pictures of $u$ and SRB measures do not
give any indication if the bifurcation from partially hyperbolic (or
Anosov) world happens. Indeed, it is very plausible that prior to
$\frac{1}{2\pi}$ no such bifurcations happen; that is, $f_{D,\eps}$ stays
Anosov with weak-strong unstable splitting for
$\eps<\frac{1}{2\pi}$. 

To provide some support we numerically calculate points of period $3$
and corresponding eigenvalues using the following procedure. We
consider a dense $2,000\times 2,000\times 2,000$  mesh of points $(x_i,y_j,z_k)$ and
apply dynamics 3 times to obtain the final point $f^3(x_i,y_j,z_k)$. If
for some $(i,j,k)$ the starting and final points end up within $D=D(x_i,y_j,z_k)<0.02$ we
adjust the coordinates $(x_i,y_j,z_k)$ coordinates to minimize the Euclidean
distance with a gradient descent method. The partial derivatives $(\frac{\partial
  D}{\partial x},\frac{\partial D}{\partial
  y},\frac{\partial D}{\partial z})$ at $(x_i,y_j,z_k)$ are
calculated numerically. Once $D<10^{-5}$, we find the cubic roots of eigenvalues of $Df^3(x_i,y_j,z_k)$.

This numerics gives 16 district eigenvalue graphs.
Indeed, Lefschetz formula yields $91$ points fixed
by $f_{*,\eps}^3$. One of these points is the fixed point $p$. The
rest give $30$ orbits of period $3$ none of which is fixed by the
involution $i$. Hence the involution breaks up these orbits into 15
pairs with identical eigenvalue data. Figure~\ref{fig45} displays
dependence of the eigenvalue data on $\eps$. We observe clear
separation of the spectrum into three bands for both conservative and
dissipative ($\eps<\frac{1}{2\pi}$) families.

\section{Background}
In this section we briefly summarize the needed background. For in depth discussions of partial hyperbolicity, SRB measures and Gibbs $u$-measures we refer the reader to~\cite{Pes, BDV, Y, PS, Dol1}.
\subsection{Anosov and partially hyperbolic diffeomorphisms}
Recall that a self-diffeomorphism $f\colon M\to M$ of a compact Riemannian manifold is called {\it Anosov} if the tangent space $T_xM$ at every $x\in M$ is split into $Df$-invariant subbundles, $T_xM=E^s(x)\oplus E^u(x)$, and $Df|_{E^s}$ is uniformly expanding while $Df|_{E^u}$ is uniformly contracting.

An important generalization is the concept of {\it partially hyperbolic diffeomorphism} $f\colon M\to M$ which assumes existence of a $Df$-invariant splitting $T_xM=E^{ss}(x)\oplus E^c(x)\oplus E^{uu}(x)$, where $Df|_{E^{ss}}$ is uniformly expanding, $Df|_{E^{uu}}$ is uniformly contracting and $Df|_{E^c}$ has intermediate growth, that is,
$$
\|Df|_{E^c}\|\cdot \|(Df|_{E^{uu}})^{-1}\|<1,
$$
and
$$
 \|Df|_{E^{ss}}\|\cdot \|(Df|_{E^c})^{-1}\|<1
$$
It is well-known that $E^{ss}$ and $E^{uu}$ integrate to foliations which we denote by $W^{ss}$ and $W^{uu}$, respectively.

For sufficiently small $\eps>0$ the diffeomorphisms $f_{*,\eps}$ given by~(\ref{D}) and~(\ref{C}) are Anosov with 2-dimensional unstable distributions. However, because $A$ has real spectrum $\lambda_1<1<\lambda_2<\lambda_3$, the unstable distribution admits a finer invariant splitting $E^c\oplus E^{uu}$ and, hence, $f_{*,\eps}$ can also be viewed as a partially hyperbolic diffeomorphism.

\subsection{SRB measures}
Informally speaking, SRB\footnote{Sinai-Ruelle-Bowen.} measures are invariant measures which are most compatible with volume when volume itself is not invariant. More precisely, consider a self-diffeomorphism $f\colon M\to M$, then an invariant measure $\mu$ is called an {\it SRB measure} (or a {\it physical measure}) if its basin of attraction has positive volume; that is, the set of points $x\in M$ such that
$$
\forall\varphi\in C^0(M)\,\,\,\lim_{n\to\infty} \frac1n\sum_{i=0}^{n-1}\varphi(f^ix)=\int_M\varphi d\mu
$$
 has positive volume. For transitive Anosov diffeomorphisms the SRB measure $\mu$ is unique and is well-understood by work of Sinai, Ruelle and Bowen. It can be characterized by the following equivalent conditions.
\begin{itemize}
\item[(C1)] $\mu$ has absolutely continuous conditionals on unstable plaques;
\item[(C2)] $\mu$ is the zero-noise limit of small random perturbations of $f$.
\end{itemize}
Our numeric calculations of the SRB measure for $f_{*,\eps}$ will rely on the second characterization which is due to L.-S. Young~\cite{Y2}. We will elaborate on it later in Section~\ref{sec_SRB}.

\subsection{Gibbs $u$-measures} The definition of Gibbs $u$-measures for partially hyperbolic diffeomorphisms comes from postulating characterization~(C1) above. Given a partially hyperbolic diffeomorphism $f\colon M\to M$, an invariant measure $\mu$ is called a {\it Gibbs $u$-measure} if it has absolutely continuous conditionals on unstable plaques. Then the density of the conditional measure on a plaque $W_f^{uu}(x,R)$ is given by
\begin{equation}
\label{eq_density}
\rho_x^{uu}(y)=\prod_{i\ge 0} \frac{Jac(f^{-1}|_{E^{uu}(f^{-i}(y))})}{Jac(f^{-1}|_{E^{uu}(f^{-i}(x))})},\, y\in W_f^{uu}(x,R)
\end{equation}

Note that in our setting the SRB measure is automatically a $u$-measure. In general, of course, the converse does not hold. Still, under additional assumptions this could be the case. For example, Bonatti and Viana showed that if $E^c$ is mostly contracting  then there are finitely many ergodic Gibbs $u$-measures which are the SRB measures~\cite{BV}. 

Dolgopyat, assuming uniqueness of the $u$-measure and that push-forwards $f_*^n\nu^{uu}$, $n\ge 0$, converge to the $u$-measure, established various limit theorems previously known in the Anosov setting~\cite{Dol2}.

\subsection{Numerics}
We have chosen the C language to implement the numerical aglorithms
for computing orbits, $u$-measures, and SRB-measures in this
study. Due to the high precision requirements in the calculation of
$u$- and SRB-measures, we employed quadmath library available in $gcc$
4.4.6. The quadruple precision \_\_float128 type provides machine
epsilon below $10^{-32}$. We relied on the Box-Muller transformation
\cite{Gaussian} of random numbers produced with a linear congruential
generator \cite{Random} to introduce Gaussian noise in the SRB
calculations. Generation of $u$-measures proved to be the most
expensive part of this numerical study but the computational cost was
fairly low, at about 1,000 CPU hours per $10^8$ points.

\section{The results}

\subsection{Numerics for the strong unstable manifold}

\label{sec_uu}

Let $f\colon\T^3\to\T^3$ be a partially hyperbolic diffeomorphism which belongs to the family~(\ref{D}) or the family~(\ref{C}) for small $\eps>0$. Here we explain numerics for the strong unstable manifold $W^{uu}_f(p)$ and present the numerical evidence supporting Conjecture~\ref{conj1} using Figures~\ref{fig7} and~\ref{fig8}. 

\setlength{\unitlength}{1mm}
\thicklines

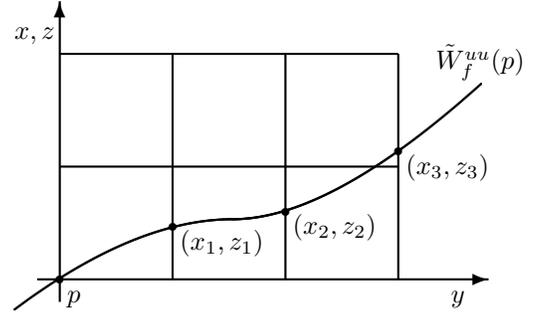
\begin{figure}

\begin{picture}(70,50)
\put(10,10){\circle*{1}}
\put(25,17){\circle*{1}}
\put(40,19){\circle*{1}}
\put(55,27){\circle*{1}}

\put(7,10){\vector(1,0){60}}
\put(10,25){\line(1,0){45}}
\put(10,40){\line(1,0){45}}
\put(25,10){\line(0,1){30}}
\put(40,10){\line(0,1){30}}
\put(55,10){\line(0,1){30}}
\put(10,7){\vector(0,1){40}}

\put(11,7){$p$}
\put(62,7){$y$}
\put(4,42){$x,z$}
\put(26,14){$(x_1,z_1)$}
\put(41,16){$(x_2,z_2)$}
\put(56,24){$(x_3,z_3)$}
\put(60,38){$\tilde W_f^{uu}(p)$}

\qbezier(4,6)(20,18)(33,18)
\qbezier(33,18)(46,18)(66,36)

\end{picture}
\caption{The lift of the strong unstable manifold and the sequence of points $\{(x_{y_0}, z_{y_0}); y_0\ge 1\}$.}
\label{fig_uu}
\end{figure}

Consider the universal cover $\R^3\simeq\{(x,y,z)\}$. Denote by $\tilde f\colon\R^3\to\R^3$ the lift of $f$ that fixes point $(0,0,0)$, which we still denote by $p$. Also denote by $\tilde W^{uu}_f(p)$ the connected component of the lift of $W^{uu}_f(p)$ which contains $p$. Then the strong unstable manifold $\tilde W^{uu}_f(p)$ can be viewed as a graph of a function $\varphi^{uu}$ defined on the $y$-axis
$$
\varphi^{uu}\colon\R\to\R^2, \,\, y\mapsto (x_y, z_y)\stackrel{\mathrm{def}}{=} \varphi^{uu}(y)
$$
as shown on Figure~\ref{fig_uu}.
For each integral $y_0$ the point $(x_{y_0}, z_{y_0})$ is the intersection point of the plane $\R^2\simeq \{y=y_0\}$ and $\tilde W^{uu}_f(p)$ and the point
$$
(x_{y_0}, z_{y_0}) \mod \Z^2
$$
is an intersection point of the 2-torus $\T^2\simeq\{y=0\}$ and $W^{uu}_f(p)$.

\begin{figure}
\includegraphics[width=70mm]{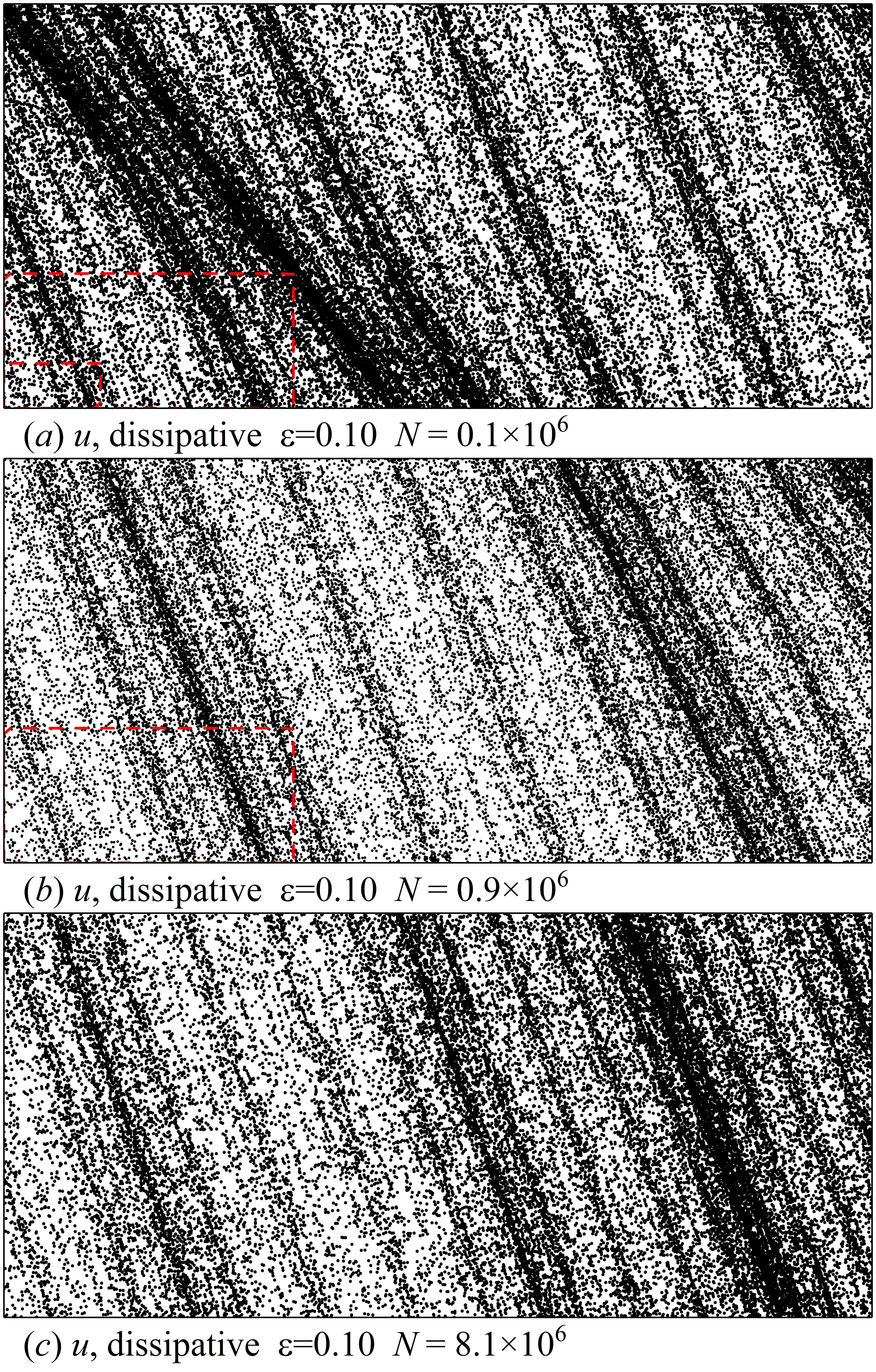}
\caption{ The intersection points for $f_D$ with $\eps=0.1$. Two lower panels are
$\times 3$ and $\times 9$ zoom-ins. The number of points is increased proportionally to the area of the domain.}
\label{fig7}
\end{figure}

To calculate $(x_{y_0}, z_{y_0})$ we carry out the following
procedure. Denote by $v_p$ the vector tangent to $\tilde W^{uu}_f(p)$
at $p$ (which is an eigenvector of $D_p\tilde f$). And let $Cv_p$ be a
vector proportional to $v_p$ whose $y$-coordinate equals to
$y_0$. Also let $\lambda_p$ be the corresponding eigenvalue,
$D_p\tilde f v_p=\lambda_p v_p$. Then to calculate the point
$q=(x_{y_0}, y_0, z_{y_0})\in \tilde W^{uu}_f(p)$ we employ the
following iterative algorithm. First we let $q^1_{-50}$ be the end
point of scaled eigenvector
$\frac{C}{\lambda_p^{50}}v_p$.\footnote{Exponent 50 is chosen so that
  in our range, $y_0\le10^8$, the distance $d(q^1_{-50}, p)$ exceeds
  the numeric error $10^{-32}$ by at least $10^8$. At the same time,
  because $\tilde W^{uu}_f(p)$ and the line spanned by $v_p$ have
  quadratic tangency, the distance from $q^1_{-50}$ to $\tilde
  W^{uu}_f(p)$ is smaller than the numeric error $10^{-32}$. } We
calculate the first approximation by using the dynamics
$$
q^1=\tilde f^{50}(q^1_{-50}).
$$

\begin{figure}
\includegraphics[width=70mm]{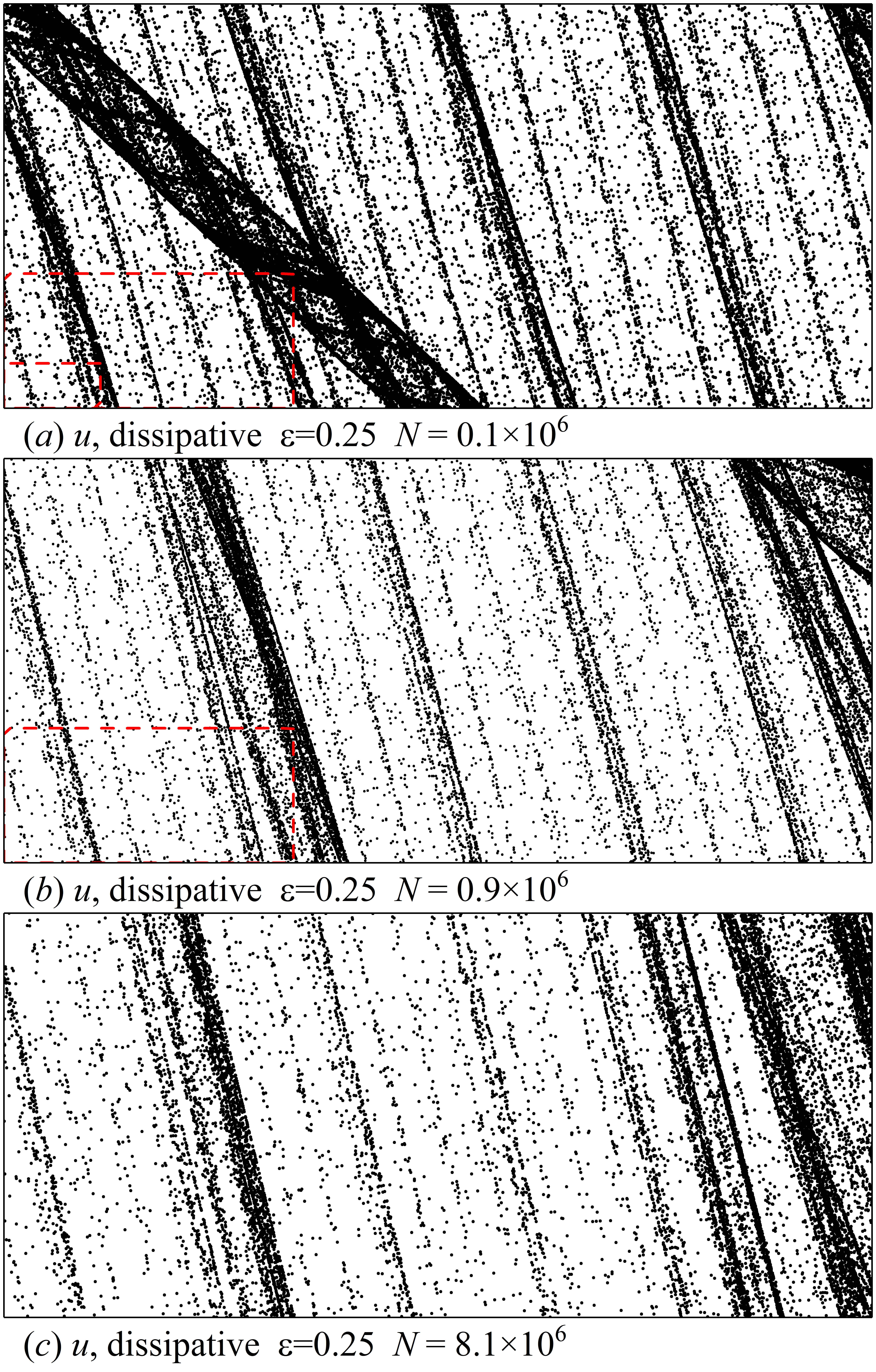}
\caption{Same as in Figure~\ref{fig7} but for $\eps=0.25$.}
\label{fig8}
\end{figure}

Then we look at the $y$-coordinate of $q^1$, compare it to $y_0$, use
a linear mixing scheme with parameter 0.1 to adjust the initial point
$q^1_{-50}$ along $v_p$ to a new initial point $q^2_{-50}$, and repeat
the procedure $n=1,000-2,000$ times to achieve convergence of $q^n$ to
the desired $q\in \tilde W^{uu}_f(p)$ within $10^{-24}$ along $\tilde W^{uu}_f(p)$.\footnote{Of
  course, exponential accumulation of error occurs during this
  procedure. The rate of error accumulation is roughly
  $\lambda_2\approx 1.55$. However note that
  $10^{-32}\times\lambda_2^{50}\approx 10^{-22}.$}

By repeating this iterative calculation we obtain the sequence of intersection points of $W^{uu}_f(p)$ and $\T^2\simeq \{y=0\}$
$$
\{(x_{y_0}, z_{y_0}); y_0\ge 1\}
$$
We calculate up to $10^8$ points in this sequence. The images we obtain clearly indicate that points cluster more for larger values of $\eps$. Still the sequence does not seem to leave any gaps in $\T^2$. This supports our density conjecture as shown on Figures~\ref{fig7} and~\ref{fig8} where as we ``zoom in" at point $p$. Zooming in does not reveal any regions free of intersection points. For the conservative family points tend to cluster much less and distribute more evenly. Hence we include the figures for the dissipative family only since they are more interesting.

\subsection{Numerics for SRB measures}
\label{sec_SRB}

We recall the description of SRB measure as zero-noise limit by L.-S. Young~\cite{Y2}. The idea is to approximate a diffeomorphism $f\colon M\to M$ by random Markov chains. To define the Markov chain  consider Borel probability measures $p(\cdot |x)$ for all $x\in M$. Given a Borel set $A\subset M$ one can think about $p(A|x)$ as the probability of sending $x$ to the set $A$. A measure $\mu$ on $M$ is {\it stationary} if
$$
\mu(A)=\int_{M}p(A |x) d\mu(x)
$$
for every Borel set $A$.

A {\it small random perturbation} of $f\colon M\to M$ is a one parameter family  of Markov chains given by transition probabilities $p^\sigma(\cdot |x)$, $x\in M$, which satisfy $p^\sigma(\cdot | x)\to\delta_{f(x)}$ as $\sigma\to 0$ uniformly in $x\in M$. (We will think of $\sigma$ as a discrete parameter.) The following properties were established in~\cite{Y2}.
\begin{enumerate}
\item If $x\mapsto p(\cdot | x)$ is continuous then a stationary measure exists;
\item If $p(\cdot | x)$ are absolutely continuous with respect to volume for all $x\in M$ then the stationary measure is also absolutely continuous;
\item If $\{p^\sigma(\cdot | x)\}$ is a small random perturbation of a diffeomorphism $f$, then all limit points of a sequence of stationary measures $\{\mu^\sigma\}$, as $\sigma\to 0$, are $f$-invariant.
\end{enumerate}

Given a measure $\nu$ on a set of diffeomorphisms $\Omega\subset\textup{Diff}(M)$ one can define the transition probabilities by
\begin{equation}
\label{eq_nu}
p(A|x)=\nu\{g: g(x)\in A\}.
\end{equation}
Now, if $\nu^\sigma\to \delta_f$ as $\sigma\to 0$ then corresponding Markov chains $\{p^\sigma(\cdot | x)\}$ yield a small random perturbation of $f$. 

Theoretical support for our computations of SRB measures which we are about to describe comes from the following theorem (which is a particular case of a more general result in~\cite{Y2}).
\begin{theorem}[\cite{Y2}]
Let $f\colon M\to M$ be a transitive Anosov diffeomorphism. There exists a $C^1$ small and $C^2$-bounded neighborhood $\Omega\ni f$ such that if $\{\nu^\sigma\}$ are Borel probability measures on $\Omega$ with $\nu^\sigma\to\delta_f$, $\sigma\to 0$, and corresponding transition probabilities $\{p^\sigma(\cdot |x)\}$ given by~(\ref{eq_nu}) are absolutely continuous then (every) sequence of stationary measures $\mu^\sigma$ converges to the SRB measure.
\end{theorem}

\begin{figure}
\includegraphics[width=70mm]{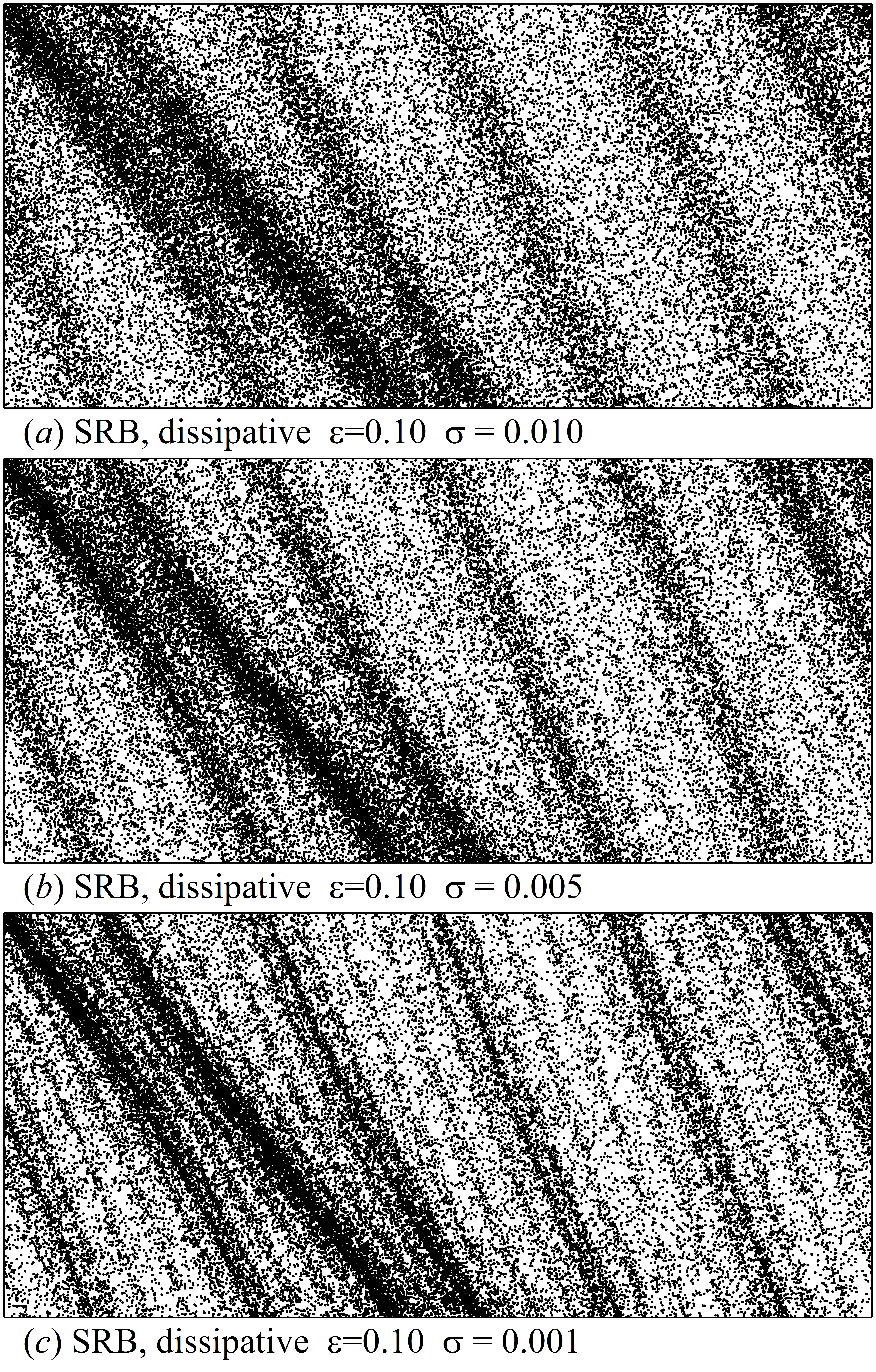}
\caption{Approximations of the SRB measure of $f_D$ with $\eps=0.1$ as $\sigma\to 0$.}
\label{fig_SRB}
\end{figure}

Hence this theorem gives a lot of credibility to numerical calculations where one applies dynamics and small random noise at each step to obtain an approximation for the SRB measure. More precisely we consider a sequence of symmetric Gaussians $\xi^\sigma$ on $\R^3$ with zero mean and standard deviation $\sigma$. Then $\xi\to\delta_{(0,0,0)}$ as $\sigma\to 0$. We define
$$
\nu^\sigma=f+(\xi^\sigma \mod \Z^3)
$$
that is, we post-compose our dynamics with a small random translation on $\T^3$. 
Then, clearly, $\nu^\sigma\to\delta_f$ as $\sigma\to 0$ and one easily sees that the transition probabilities are absolutely continuous. Hence the theorem above applies. (Technically, we also need to truncate Gaussians to ensure $C^1$-smallness of the perturbation, but practically this makes no difference as we are interested in very small $\sigma$.)

The numeric scheme is as follows. We begin with a random point $q_0$ on $\T^3$ and generate a $\sigma$-approximation of the SRB measure by consecutive application of $f$ and addition of Gaussian noise $\xi^\sigma$. That is,
$$
q_{i+1}=f(q_i)+\xi^\sigma
$$
Note that we only work in the parameter range $\sigma\gg 10^{-32}$ so that the numeric error in calculation of $f$ is much smaller than the (small) random noise. Therefore, exponential accumulation of the numeric error is not of any concern. On Figure~\ref{fig_SRB} we display several approximations for different values of $\sigma$. For all further SRB measures numerics, which we need for comparisons with $u$-measures numerics, we use $\sigma=10^{-29}$.

\subsection{Comparing Gibbs $u$ and SRB measures}

Consider the averaged Dirac measures
\begin{equation}
\label{eq_dirac_u}
\Sigma^u=\frac1N\sum_{y_0=1}^N \delta_{(x_{y_0}, z_{y_0})}
\end{equation}
corresponding to  the $u$-measure, and
\begin{equation}
\label{eq_dirac_SRB}
\Sigma^{\textup{SRB}}=\frac1N\sum_{q_i\in S} \delta_{q_i}
\end{equation}
corresponding to the SRB measure, where $S$ is the slice $\{(x,y,z\in\T^3: \, -0.005\le y\le 0.005\}$ which contains $N$ points. 
For the dissipative family point distributions $\Sigma^u$ and $\Sigma^{\textup{SRB}}$ visually coincide for all parameters in our range $\eps\in[0, 0.25]$ (see  Figure~\ref{fig3}). On the other hand, for the conservative family, the SRB measure is the uniform Lebesgue measure as it supposed to be, while the $u$-measure appears to be an absolutely continuous measure with a non-constant density, as one can see on the top panel of Figure~\ref{figw1}. The ``non-uniformity" increases as we increase $\eps$.  The explanation for this discrepancy is that~(\ref{eq_dirac_SRB}) gives the (approximation of)  true conditional measure on $\T^2$ of the SRB measure, while~(\ref{eq_dirac_u}) does not give (an approximation of) the conditional of $f_*^n\nu^{uu}$. Hence we proceed with the numeric calculation of the true conditional of $f_*^n\nu^{uu}$ on $\T^2$ and present the numeric evidence that the measures indeed coincide.
\begin{remark}
Note however that the conditional of $f_*^n\nu^{uu}$ on $\T^2$ is absolutely continuous with respect to~(\ref{eq_dirac_u}). Hence, if $\Sigma^u$ converges to an absolutely continuous measure (which we numerically verified by using histograms)  then $f_*^n\nu^{uu}$ converges to an absolutely continuous measure on $\T^3$ as $n\to\infty$. And, since this measure is invariant, it must be the volume. In view of this remark our further numeric verification of convergence of $f_*^n\nu^{uu}$  to volume becomes somewhat redundant. However we still find it important to have direct numeric evidence.
\end{remark}
\begin{remark}
By analyzing distribution functions of the Dirac averages~(\ref{eq_dirac_u}) and~(\ref{eq_dirac_SRB}) in the dissipative family we can also very clearly conclude that $\Sigma^u$ and  $\Sigma^{\textup{SRB}}$ do not converge to the same measure on $\T^2$. Hence, as to be expected, the above discussion also applies to the dissipative family. However visually the point distributions $\Sigma^u$ and  $\Sigma^{\textup{SRB}}$ are identical in this case. This happens because for singular measures, when looking at the pictures of approximating point distributions we can only see the measure class rather than the measure itself.
\end{remark}

\subsection{The conditional measure on  $\T^2\simeq\{y=0\}$ for Gibbs $u$-measure}

Now we explain precisely our numerics for the conditional of the Gibbs $u$-measure. Let $r\in \tilde W_f^{uu}(p)$ be a point very close to $p$. Consider the Lebesgue measure on $\tilde W_f^{uu}(p)$ induced by the canonical flat Riemannian metric on $\mathbb R^3$. And denote by $\nu^{uu}$ the normalized Lebesgue measure supported on the strong unstable plaque $[p,r]^{uu}\subset\tilde W_f^{uu}(p)$. Then, by using calculus, the density of $f_*^n\nu^{uu}$ with respect to the Lebesgue measure on $\tilde W_f^{uu}(p)$ is given by
$$
\rho(q)=Jac(f^{-n}|_{E^{uu}(q)}), \,\, q\in [p, f^n(r)]^{uu}
$$
This Jacobian density can be easily evaluated numerically because, as we explained in Section~\ref{sec_uu}, we can accurately calculate points $q$ on $\tilde W_f^{uu}(p)$ together with their preimages under $f^{-n}$, $n\le 50$. Hence to find $\rho(q)$ approximately we look at points $q-\Delta q, q+\Delta q$ on $\tilde W_f^{uu}(p)$ and their preimages; and then evaluate the Jacobian numerically in the obvious way.\footnote{We use $\Delta q=10^{-10}$.}

Recall that we need to further take the conditional measure of $f_*^n\nu^{uu}$ on  $\T^2\simeq\{y=0\}$. Then, one can easily see (for example, by taking the limit as the width of the $\T^2$-slice goes to zero) that the expression for the conditional measure at the intersection point $(x_{y_0}, y_0, z_{y_0})$, $y_0\in\Z_+$, depends on the angle between $\tilde W_f^{uu}(p)$ and $\T^2$ at $(x_{y_0}, y_0, z_{y_0})$. Namely, one has the following formula for the conditional of $f_*^n\nu^{uu}$ on  $\T^2\simeq\{y=0\}$ 
$$
\Sigma^u_{\rho a}=\frac{1}{W(N)}\sum_{y_0=1}^N \rho(x_{y_0}, y_0, z_{y_0}) a(x_{y_0}, z_{y_0})\delta_{(x_{y_0}, z_{y_0})},
$$
where $\rho$ was defined above, 
$$
W(N)=\sum_{y_0=1}^N \rho(x_{y_0}, z_{y_0}) a(x_{y_0}, y_0, z_{y_0});
$$
and the ``angle weight" is defined by
$$
a(x_{y_0}, z_{y_0})=\frac{1}{\langle v^\perp, v^{uu}\rangle},
$$
where $v^\perp$ is unit vector at $(x_{y_0}, y_0, z_{y_0})$ perpendicular to $\T^2$ and $v^{uu}$ is the unit vector at $(x_{y_0}, y_0, z_{y_0})$ tangent to  $\tilde W_f^{uu}(p)$. Again, coefficient $a$ is easy to calculate since we can numerically calculate the tangent vectors $v^{uu}$.

The density weights $\rho$ and $a$ have different genesis. Hence, for the purpose of analyzing $\Sigma^u_{\rho a}$ we also introduce the ``component" Dirac averages
$$
\Sigma^u_{\rho}=\frac{1}{W}\sum_{y_0=1}^N \rho(x_{y_0}, y_0, z_{y_0})\delta_{(x_{y_0}, z_{y_0})},
$$
and
$$
\Sigma^u_{a}=\frac{1}{W}\sum_{y_0=1}^N  a(x_{y_0}, z_{y_0})\delta_{(x_{y_0}, z_{y_0})},
$$
which are normalized by the corresponding total weight $W$. 

\subsection{Comparing Gibbs $u$ and SRB measures numerically}

\begin{figure}
\includegraphics[width=79mm]{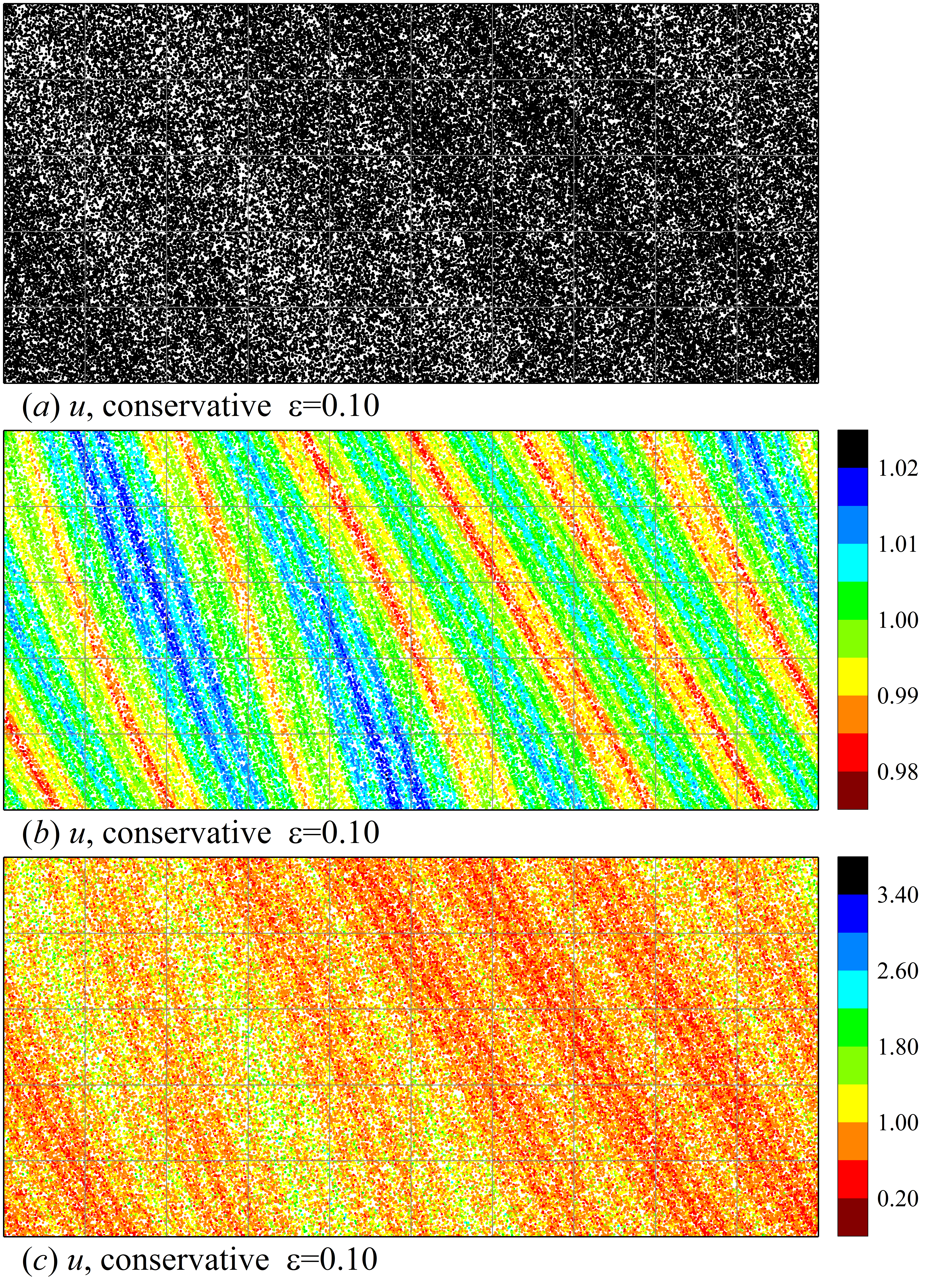}
\caption{Dirac averages $\Sigma^u$, $\Sigma^u_a$ and $\Sigma^u_\rho$, $N=200,000$. Note that, for example, in right lower corner point distribution $\Sigma^u$ takes large values,  $a$ is above average, and $\rho$ takes small values.}
\label{figw1}
\end{figure}

\begin{figure}
\includegraphics[width=85mm]{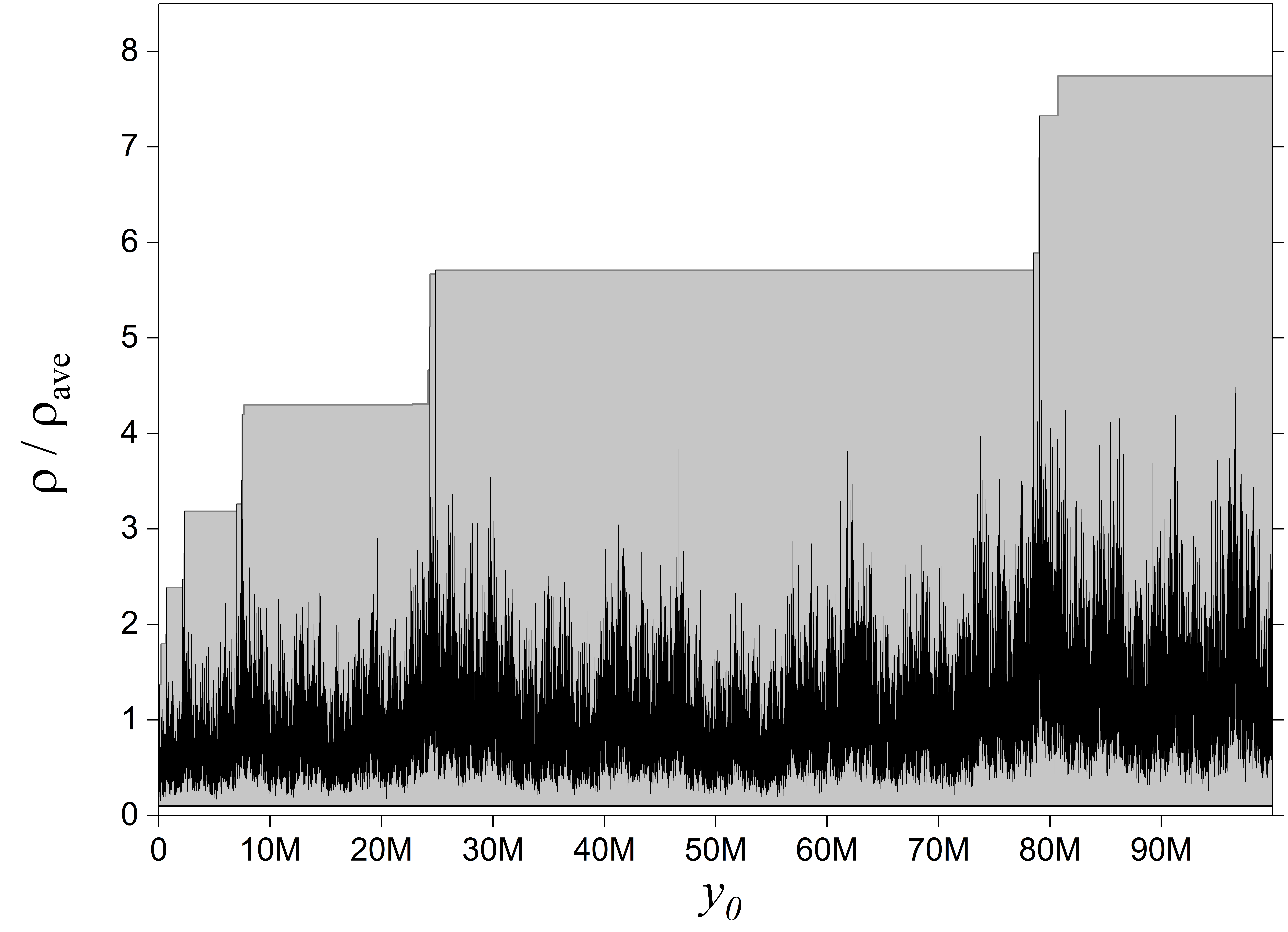}
\caption{Dependence of $\rho/\rho_{ave}$ on $y_0$ for $u$-measure of
  $f_{C,\eps=0.1}$. $\rho_{ave}$ is calculated over the full range of
  $N=10^8$ points. The minimum $0.1245$ is achieved at 1.}
\label{figw2}
\end{figure}

\begin{figure}
\includegraphics[width=85mm]{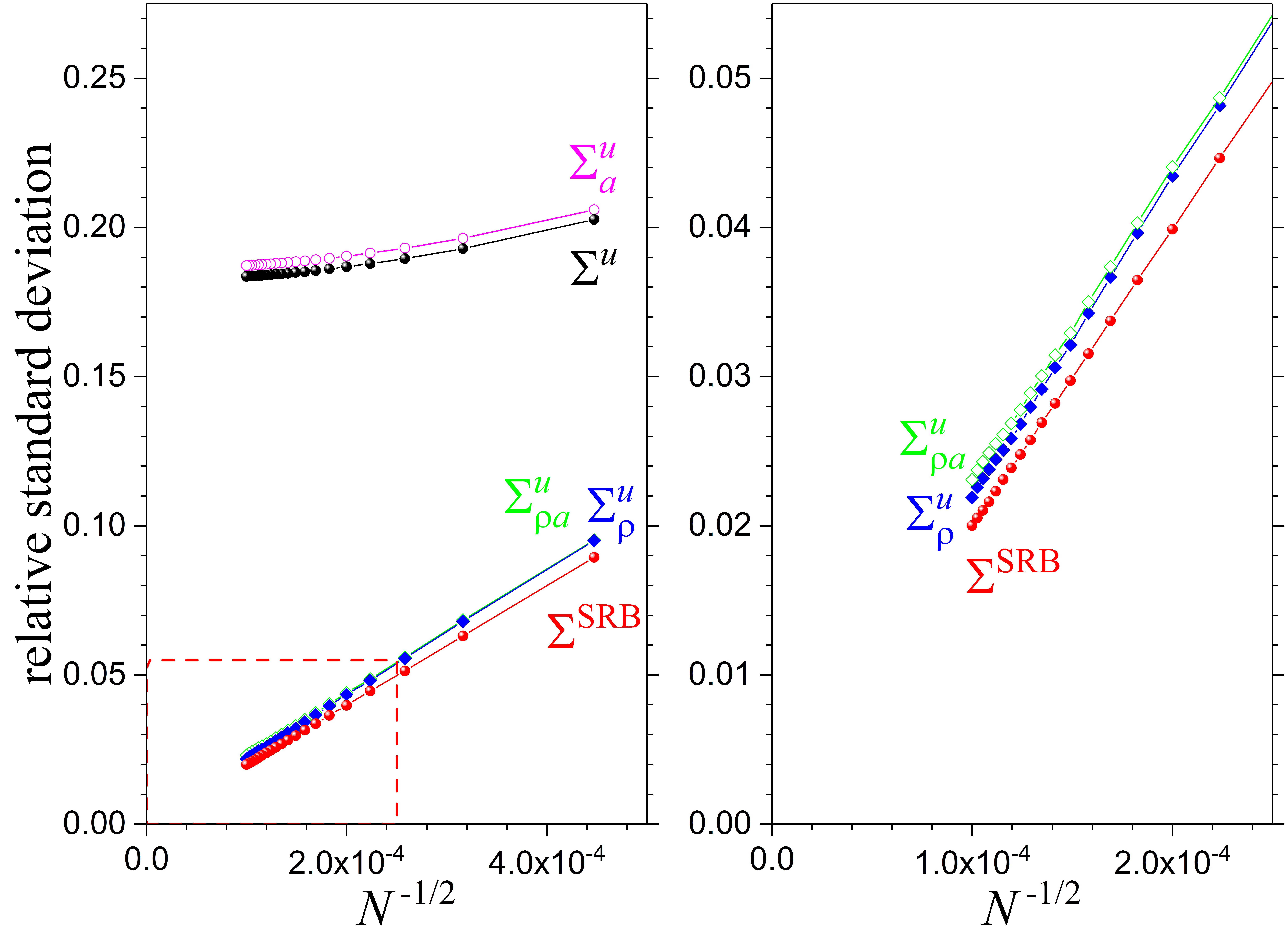}
\caption{Relative standard deviations for $f_{C,\eps=0.1}$ point
  distributions calculated for $N$ between $5$ to $100$ million. The
  zoom-in on the right shows that both SRB measure and $u$ measure
  (when properly weighted) converge to uniform distributions
  approximately as $N^{-1/2}$. }
\label{figw3}
\end{figure}

Our calculations of $\rho$, $a$ and point distributions $\Sigma^u$, $\Sigma^u_\rho$, $\Sigma^u_a$ and $\Sigma^u_{\rho a}$ are summarized on Figures~\ref{figw1},~\ref{figw2} and~\ref{figw3}.

\captionsetup{width=79mm}

On Figure~\ref{figw1} the values of weights are coded in color. The average value is normalized to equal 1. By the definition the ``angle weight" $a=a(x_{y_0}, z_{y_0})$ is a continuous function on $\T^2$. Notice that $a$ (middle panel) varies only slightly, within 2\% of the average. On the other hand, on the bottom panel shows that $\rho$ varies a lot. Further, the graph of $\rho=\rho(x_{y_0}, y_0, z_{y_0})=\rho(y_0)$ given on Figure~\ref{figw2} shows that $\rho$ is unbounded (that is, if normalize $\rho$ so that $\rho(0)=1$ then $\rho$ is unbounded function of $y_0$) and have certain ``self-similar" structure. 

By examining Figure~\ref{figw1} one can see that $\rho$ smoothes out the point distribution $\Sigma^u$, that is, it makes $\Sigma^u_\rho$ more uniform than $\Sigma^u$. Curiously, and we have no good explanation for this, the ``angle weight" $a$ makes point distribution less uniform, but, as we remarked before, $a$ has a very small effect on the distribution.

\captionsetup{width=85mm}

In order to quantify these observations coming from Figure~\ref{figw1} we use a $200\times 200$ square grid to partition $\T^2$ into $40,000$ bins. For each bin $\mathcal B$ we calculate its total weight
$$
w^u(\mathcal B)=\#\{y_0\le N: (x_{y_0}, z_{y_0})\in\mathcal B\}
$$
as well as the weight adjusted by $\rho$
$$
w^u_\rho(\mathcal B)=\sum_{\substack{y_0=1,\\ (x_{y_0},z_{y_0})\in\mathcal B}}^N \rho(y_0).
$$
And weights $w^{\textup{SRB}}(\mathcal B)$, $w^u_a(\mathcal B)$ and $w^u_{\rho}(\mathcal B)$ are defined analogously. Further we calculate relative standard deviation in order to have a single number which measures closeness to the uniform distribution
$$
RSD^u=\frac{1}{\bar w^u}\left(\frac{1}{40,000}\sum_{\mathcal B}(w^u(\mathcal B)-\bar w^u)^2\right)^{\frac12},
$$
where $\bar w^u$ is the average of the weights. Analogously we have relative standard deviations $RSD^{\textup{SRB}}$, $RSD^u_\rho$, $RSD^u_a$ and $RSD^u_{\rho a}$. The dependence of relative standard deviations on the number of points in the range $N=10^6,\ldots 10^8$ is shown on Figure~\ref{figw3}. Indeed, we see that $RSD^\textup{SRB}$, $RSD^u_\rho$ and $RSD^u_{\rho a}$ decay to zero roughly proportionally to $\frac{1}{\sqrt N}$. Unfortunately we cannot differentiate between $RSD^u_\rho$ and $RSD^u_{\rho a}$.

If we denote by $F^u_*\colon [0,1]^2\to [0,1]$ the distribution function of $\Sigma^u_*$, $*=\rho, \rho a$, given by
$$
F^u_*(c,d)=\Sigma^u_*([0,c]\times[0,d])
$$
then weak$^*$ convergence of $\Sigma^u_*$ to the Lebesgue measure is equivalent  to convergence of $F^u_*(c,d)$ to $cd$ for all $(c,d)\in[0,1]^2$ as $N\to\infty$. We remark that convergence of $RSD^u_*$ to $0$ is equivalent to
$$
F^u_*(c,d)\to cd,\,\,\, (c,d)=\left(\frac{i}{200}, \frac{j}{200}\right), \,i,j=0,\ldots 200.
$$
Such convergence of distribution functions is known in statistics as Kolmogorov-Smirnov test~\cite{Smi}.

Finally let us mention that we have also performed similar numerics, such as the Kolmogorov-Smirnov test, for the dissipative family and the results are similar. It is more difficult to compare $\Sigma^\textup{SRB}$ and $\Sigma^u_{\rho a}$ in this case because the SRB measure is not Lebesgue. The difficulty comes from the fact that $\Sigma^\textup{SRB}$ is defined by~(\ref{eq_dirac_SRB}) using the slice $S$ of thickness $\Delta y=10^{-2}$. In conservative case the value of $\Delta y$ is irrelevant, but in the dissipative case the restriction of the SRB measure to the slice is no longer a product measure. Hence we also need to let $\Delta y\to 0$ in order to approximate the conditional measure on $\T^2$. The additional parameter $\Delta y$ makes numerics even more involved and we did not fully pursue it.

{\it Acknowledgments.} 
A.G. would like to thank Yakov Pesin who introduced him to questions in the spirit of our Conjecture~\ref{conj2} in his 2004 dynamics course. Also A.G. would like to thank Aleksey Gogolev who performed  initial numerical experiments and created the first set of beautiful pictures back in 2008. During final stages of preparation of this paper discussions with Federico Rodriguez Hertz were very useful. We also would like to acknowledge helpful feedback from Dmitry Dolgopyat, Yi Shi and Rafael Potrie.

A.G. was partially supported by NSF grant DMS-1204943. I.M. and A.N.K. gratefully acknowledge NSF support (Award No. DMR-1410514).

\medskip
\medskip
\medskip
Andrey Gogolev, \\
Mathematical Sciences,\\
SUNY Binghamton, N.Y., 13902, USA\\
{\tt agogolev@math.binghamton.edu }\\
\\

\noindent Aleksey Kolmogorov, \\
Physics Department,\\
SUNY Binghamton, N.Y., 13902, USA\\
{\tt kolmogorov@binghamton.edu }\\
\\

\noindent Itai Maimon\\
Mathematics and Physics\\
SUNY Binghamton, N.Y., 13902, USA\\
{\tt imaimon1@binghamton.edu }

\end{document}